\documentclass[12pt,a4paper]{article}

\usepackage[french,english]{babel}
\usepackage{amsfonts}
\usepackage{amssymb}
\usepackage{amsthm}
\usepackage[T1]{fontenc}
\usepackage[centertags]{amsmath}
\usepackage{newlfont}

\usepackage{epsfig}
\usepackage{float}

\newtheorem{defi}{Definition}
\newtheorem{defiprop}{Definition-Proposition}
\newtheorem{rem}{Remark}
\newtheorem{prop}{Proposition}
\newtheorem{theo}{Theorem}

\newtheorem{coro}{Corollary}

\newcommand{\set}[1]{\left\{#1\right\}}
\newcommand{\A}{\mathbb A}
\newcommand{\al}{\alpha}
\newcommand{\be}{\beta}
\newcommand{\g}{\gamma}

\begin{document}

\title{COMBINATORIAL TANGENT SPACE AND RATIONAL SMOOTHNESS OF SCHUBERT VARIETIES}

\author{Stéphane Gaussent}

\date{Laboratoire Géométrie, Topologie et Algèbre\\ Université Montpellier
2, Case 051, Pl. E. Bataillon\\ 34095 MONTPELLIER Cédex 05 \\
courriel : gaussent@math.univ-montp2.fr}

\maketitle

\selectlanguage{french}
\begin{abstract}

A la suite de Contou Carrere [CC], nous considérons la résolution
de Bott-Samelson d'une variété de Schubert comme une variété de
galeries de l'immeuble de Tits associé à la situation. Nous
montrons que la lissité rationnelle d'une telle variété est codée
par un sous-espace de l'espace tangent de Zariski appelé, l'espace
tangent combinatoire. Nous utilisons pour cela une
caractérisation de la lissité rationelle d'une variété de
Schubert introduite par Carrell et Peterson [CP].

\end{abstract}
\selectlanguage{english} \abstract{Following Contou Carrere [CC],
we consider the Bott-Samelson resolution of a Schubert variety as
a variety of galleries in the Tits building associated to the
situation. We prove that the rational smoothness of a Schubert
variety can be expressed in terms of a subspace of the Zariski
tangent space called, the combinatorial tangent space. For this,
we use a characterization of rational smoothness of a Schubert
variety introduced by Carrell and Peterson [CP].}



\newpage

\section{Introduction}

In this paper, we introduce a combinatorial tool in the study of
the singular locus of a Schubert variety, called the
combinatorial tangent space. This space is built from a
Bott-Samelson resolution of such a variety and is extracted from
the work of Contou Carrere [CC]. Using a combinatorial result of
Carrell and Peterson [CP, Theorem C] which characterizes, via the
Kashdan-Lusztig polynomials, the rational smoothness of a
Schubert variety, we show our main result (Section 4, Theorem 4):

\bigskip
{\bf Theorem.}  {\it The Schubert variety $\overline\Sigma (B,
w)$ is rationally smooth at a $T-$fixed point $x =  u(P)$ if and
only if for all $u\leq v\leq w$, the dimension of the
combinatorial tangent space of $\overline\Sigma (B,w)$ at $v(P)$
is equal to the dimension of the variety.}

\bigskip

Here is an outline of the paper, in $\S$2, we state some
definitions and some combinatorial settings. In particular, we
recall from [CC], the universal Schubert scheme and the
construction of the Bott-Samelson variety as a variety of
galleries in the Tits building associated to any semi-simple
(adjoint) group $G$. In $\S$3, we give the definition of the
combinatorial tangent space and state some of these properties.
In $\S$4, we relate explicitly the combinatorial tangent space
and the rational smoothness of a Schubert variety. And we show
how this space is connected to the set of $T-$invariant curves in
such a variety.

\bigskip
Althought the construction of the combinatorial tangent space of a
Schubert variety is presented here in the theory of semi-simple
groups, this construction is still valid in the more general
setting of Kac-Moody groups and their Schubert varieties [G].
Indeed, in this context, Bott-Samelson and Schubert varieties can
be defined [Ma] and they have a combinatorial interpretation, of
course, by taking into account the building associated to such a
group. The notion of "épinglage" (cf. Section 2 below) does not
exist for a Kac-Moody group $G$, but we still have a basis
(infinite this time) of the tangent space of $G/B$ at any
$T-$fixed point (given by the presentation of the Kac-Moody
algebra). Moreover, the characterization of rational smoothness
we use generalizes to this case ([KL] and [C1]). So does our main
result : the rational smoothness of the Schubert varieties is
characterized by the combinatorial tangent space.

This approach, i.e. studying the smoothness of a Schubert variety
via the Bott-Samelson resolution (viewed as variety of galleries),
follows the work of Contou Carrere [CC]. And it has many others
applications, some of them could be found in my Ph. D. thesis [G]
that I am working out at Université Monpellier 2 (France).

\bigskip
I would like to thank my supervisor, C. Contou Carrere, for his
guidance through these topics. I also thank Ph. Elbaz-Vincent and
P. L. Montagard for their comments and their careful readings of
the manuscript and Ph. Malbos for his initiation into Pstricks.


\section{Settings and Definitions}

In this section, we use partially the notations of [SGA3] and we
present the Schubert variety, following Contou Carrere [CC], in a
non usual way, which is quite near from the point of view of
Kazhdan and Lusztig [KL].

Let $G=(G,T,M,R)$ be a semi-simple (adjoint) split group over an
algebraically closed field $k$. That is, $T$ is a maximal torus of
$G$, $M$ is an abelian group such that $Hom_{k-gr} \big (T,\mathbb
G_m(k)\big )\simeq M_k$ and $R$ is a root system of $(G,T)$. We
have the following decomposition of ${\mathcal G}=Lie (G)$ under
the $T-$adjoint action :

$${\cal G} = {\cal T} {\mathbf\oplus} \bigoplus_{\alpha\in R} {\cal G}^{\alpha},$$
where the ${\cal G}^{\alpha}$'s are one dimensional vector spaces.

Let $\varepsilon = (R_0, \{ X_{\alpha} \}_{\alpha\in R_0})$ be an
``épinglage'' of $G$. That is, we choose a set $R_0$ of simple
roots in  $R$ (this choice corresponds to the choice of a Borel
subgroup $B$ of $G$ such that $T\subset B\subset G$) and for each
$\alpha\in R_0$, we choose $X_{\alpha}$ a generator of the vector
space ${\cal G} ^{\alpha}$.

Let $W= Norm_G(T)/T$ denote the Weyl group of $(G,T)$ and let $S$
be a set of generators of $W$ adapted to the choice of $B\subset
G$ (in particular, there is a bijection between $R_0$ and $S$, we
thus index the elements of $S$ by simple roots). Moreover, for a
part $t$ of $S$, $W_t\subset W$ is the subgroup of $W$ generated
by $\{s_\alpha\}_{\alpha\in t}$.

\subsection{Universal Schubert Scheme}

Following [CC, Part I, $\S$4], we introduce the universal
Schubert scheme, we refer to this work for more details and a
full description of these notions.

For any set $E$, $\mathcal P(E)$ stands for the set of all his
parts. Let $Par(G)$ denote the $k-$variety of all the parabolic
subgroups of $G$. Herewith the ``épinglage'' of $G$, we get a
decomposition of this variety:

$$Par (G) = \coprod_{t\in {\cal P}(R_0)} Par_t(G),$$ where
$Par_t(G)$ is defined as the variety of parabolic subgroups of
type $t$. The variety of Borel subgroups of $G$ is then $Bor(G) =
Par_{\emptyset}(G)$.

Let $\Sigma$ denote the variety $Par(G)\times_k Par(G)$. The group
$G$ acts diagonally on $\Sigma$ and the quotient is denoted by
$t.st = \Sigma / G$. This is the constant $k-$scheme :

$$\coprod_{t,t'\in {\cal P}(R_0)} W_t\backslash W / W_{t'}.$$

In the sequel of what we wrote, we can define an immersion :

$$\Sigma\hookrightarrow Par(G)\times_k t.st\times_k Par(G),$$
and state the following definition.

\begin{defi}
We call $\Sigma$ the universal Bruhat cell and $\overline\Sigma$,
the scheme theoretic image of $\Sigma$ in $Par(G)\times_k
t.st\times_k Par(G)$ by this immersion, is called the universal
Schubert scheme.
\end{defi}

\begin{rem}
If we fix a parabolic subgroup $P$ of type $t_P$ such that
$P\supset B$, where $B$ is the Borel subgroup of $G$ given by the
``épinglage'', and an element $\overline w\in W/W_{t_P}$, then by
taking fibres in the previous definition, we recover the Bruhat
cell $\Sigma (B,\overline w)$ and the Schubert variety
$\overline\Sigma (B,\overline w)$ :

$$\overline\Sigma (B,\overline w) \hookrightarrow \{ B
\}\times\{\overline w \}\times Par_{t_P}(G) \simeq G/P.$$

The cell $\Sigma (B,\overline w)$ is the variety of parabolic
subgroups of type $t_P$ of $G$ in relative position $\overline w$
with $B$ and $\overline\Sigma (B,\overline w)$ is the variety of
parabolic subgroups of type $t_P$ of $G$ in relative position
$\overline u$ with $B$ for $\overline u \leq \overline w$ in
$W/W_{t_P}$ (here $\leq$ represent the relative Chevalley-Bruhat
order).
\end{rem}

In the following, we will only consider Schubert varieties in
$Par_{t_P}(G)\simeq G/P$, where $P\supset B$ is a parabolic
subgroup of $G$ of type $t_P$.


\subsection{Tits Building}

Now, let us go into a succint description of a combinatorial
object naturally associated to a semi-simple group, the Tits
building (actually, this construction is valid without any
hypothesis on the ground field $k$).

We denote by $\Delta (G)$ the set of all parabolic subgroups of
$G$ ordered by the opposite relation of the inclusion between
parabolic subgroups. This operation endowes this set with a
structure of simplicial complex. The variety $Par(G)$ and the
complex $\Delta (G)$ have the same underlying set, but different
structures.

As we have fixed a maximal torus, $T\subset G$, we can define a
subcomplex $\Bbb A$ associated to $T$, called an apartment, as
follows :

$$\Bbb A = \{ P\in Par (G),\ T\subset P\},$$ if $P\in\Bbb A$ then
there exists a parabolic subgroup $Q$ containing $B$ and $n\in N=
Norm_G (T)$ such that $P=nQn^{-1}$.

And the building can be expressed as $\Delta (G) = G\times \Bbb A
/\sim \,$ where $\sim$ is the equivalence relation :

$$(x,P)\sim (x',P') \Leftrightarrow
\left\{
\begin{array}{l}
P=P'\\
x^{-1}x'\in P.
\end{array}
\right. $$ For each maximal torus in $G$ there is an apartment
corresponding to it and in order to recover $\Delta(G)$ we paste
all these subcomplexes with respect to the equivalence relation
$\sim$. This description of $\Delta(G)$ as a ``union'' of
apartments provides a better handing to deal with the building.

\bigskip

Elements of $\Delta (G)$ are called faces, to each parabolic
subgroup $P$ of $G$ is associated a face, denoted by $F_P$. A
chamber of $\Delta (G)$ is a maximal face. The chambers are
associated to the Borel subgroups of $G$.

The group $W$ acts on the faces of $\Bbb A$ and in a simply
transitively way on the set of all the chambers of $\Bbb A$.
Moreover, we can define a distance between two chambers of $\Bbb
A$ :

$$d(C,C') = l(w_{CC'}),$$ where $w_{CC'}$ is the element of $W$
which maps $C$ to $C'$ and $l(w)$, for $w\in W$ is the length of
any reduced decomposition of $w$.

To each face $F_P$ contained in the chamber $F_B$, we associate
its type :

$$typ(F_P) = t\subset S,$$ where $t$ is defined by $W_t=
Stab_W(F_P)$. Therefore, we get $typ (F_B)= \emptyset$ since
$Stab_W (F_B)=\{ 1\}=W_{\emptyset}$. The action of $W$ allows us
to define the type of each face of $\Bbb A$. Moreover, thanks to
the retraction of $\Delta (G)$ onto $\Bbb A$ with centre $F_B$ [T,
3,3, p.42], we can define the type of all the faces (as there is a
bijection between $R_0$ and $S$, the two notions of type
mentionned are the same).

The set of types, denoted by $typ\big (\Delta (G)\big )$ is still
a complex and $typ : \Delta (G) \to typ\big (\Delta (G)\big )$ is
a morphism of complex.

\bigskip

A {\it folding} of $\Bbb A$ is an idempotent and type preserving
morphism of complex $\phi : \Bbb A \to \Bbb A$ such that each
chamber $C$ belonging to $\phi (\Bbb A)$ is the image of exactly
two chambers of which itself [T, 1,8, p.7]. The image of a folding
$\phi (\Bbb A)$ is called a {\it root} of $\Bbb A$.

If $\alpha$ is a {\it root} of $\Bbb A$, we define $M_{\alpha}$,
the {\it wall associated to $\alpha$}, as the subcomplex of $\Bbb
A$ composed by the faces $F$ such that there are two adjacent
chambers $C$ and $C'$ (i.e. $C\cap C'$ is a codimension 1 face of
$C$ and of $C'$) with $C\in \alpha, C'\not\in\alpha$ and $F\subset
C\cap C'$.

If we look at a graphic representation of $\Bbb A$ as the root
system of $(G,T)$ (for instance see [BOU, Planche X]), the {\it
roots} correspond to the roots of $(G,T)$ and the walls correspond
to the hyperplans related to the reflections. To each wall
$M_{\beta}$ in $\Bbb A$, we can associate two foldings
$\phi_{\beta}$ and $\phi_{-\beta}$, corresponding to the roots
$\beta$ and $-\beta$.

Moreover, we say that a folding $\phi$ is {\it towards a chamber
$C$} if $C$ belongs to the image of $\phi$. For example, the
foldings towards $F_B$ are indexed by the positive roots.


\subsection{Generalized Galleries}

The notion of gallery that we are using is a bit more general than
the one presented in Tits' book [T]. So, in order to make a
difference, we will use the expression $g-$gallery to denote such
a gallery.

\begin{defi}[{[CC, Part I, $\S$2]}]

Let $(E,<)$ be an ordered set. A sequel $(e_n,...,e_0)$ of
elements of $E$ is a $g-$gallery if those elements fit into one of
the following situations :

$$
\begin{array}{c}
(e_n>e_{n-1}<\cdots >e_1<e_0)\\
(e_n<e_{n-1}>\cdots >e_1<e_0)\\
(e_n>e_{n-1}<\cdots <e_1>e_0)\\
(e_n<e_{n-1}>\cdots <e_1>e_0).
\end{array}
$$
\end{defi}

For example, let $g$ be a $g-$gallery of $\Delta (G)$ of which
faces verify the following incidence relations :

$$g = (F_r \supset F'_r \subset F_{r-1} \supset \cdots \supset
F'_1 \subset F_0 \supset F'_0).$$ We define the type of $g$ to be
the $g-$gallery of types $\tau$ :

$$\tau = typ(g)=(t_r \subset t'_r \supset t_{r-1} \subset \cdots
\subset t'_1 \supset t_0 \subset t'_0),$$ where each $t_i$ (resp.
$t_j'$) is the type of the corresponding face.

Actually, $g$ corresponds to a configuration of parabolic
subgroups verifying the following inclusions :

$$g =
(Q_r \subset P_r \supset Q_{r-1} \subset \cdots \subset P_1
\supset Q_0 \subset P_0),$$ where the $Q_j$'s are of type $t_j$
and the $P_i$'s, of type $t'_i$.

\begin{defi}
If $F_r$ and $F'_0$ are two faces of $\Bbb A$ (resp. of $\Delta
(G)$), we denote by $$\Gamma (\Bbb A ,\tau ,F_r,F'_0)\quad
({resp.}\ \Gamma (\Delta (G),\tau ,F_r,F'_0))$$ the finite set
(resp. the set) of all $g-$galleries contained in $\Bbb A$ (resp.
in $\Delta (G)$), of type $\tau$, of source $F_r$ and of target
$F'_0$.
\end{defi}

\bigskip
Now, we present the notion of minimality for $g-$galleries.

For two faces $F$ and $F'$ of $\Bbb A$ such that $F'\subset F$, we
denote by ${\cal M}_{F'}(F)$ the set of walls $M$ such that $F'\in
M, F\not\in M$.

\bigskip
\begin{defi}[{[CC, Part I, $\S$5]}]
A g-gallery $\gamma = (F_r \supset F'_r \subset F_{r-1} \supset
\cdots \supset F'_1 \subset F_0 \supset F'_0)\in \Gamma (\Bbb A
,\tau ,F_r,F'_0),$ between the chamber $F_r=F_B$ and the face
$F_0'$ is said to be minimal if :

i) for $r\geq i\geq 1$ the sets ${\cal M}_{F'_i}(F_{i-1})$ are
disjoint ;

ii) ${\cal M}(F_B,F'_0)= \cup_{r\geq i\geq 1} {\cal
M}_{F'_i}(F_{i-1})$, where ${\cal M}(F_B,F'_0)$ denote the set of
walls which separate $F_B$ from $F'_0$.
\end{defi}

For a discussion on further properties of minimal $g-$galleries,
we refer to [CC].


\subsection{Chamber Galleries}

\begin{defi}
A g-gallery $g = (F_r \supset F'_r \subset F_{r-1} \supset \cdots
\supset F'_1 \subset F_0 \supset F'_0)$ is called a chamber
g-gallery if

i) for $r\geq i\geq 0$, $F_i$ is a chamber ;

ii) for $r\geq j\geq 1$, $F'_j$ is a face of codimension 1 in
$F_j$ and $F_{j-1}.$

\end{defi}

The chamber $g-$galleries are nearly galleries in the sense of
Tits [T], except that we allow their target face to be of any type
(instead of being a chamber). By the way, we will make the abuse
of language and call our $g-$galleries, simply galleries.

So, the minimality of a chamber gallery $g$ is expressed in the
same way as in [T]. Keeping the notations as in definition 5, $g$
is minimal if :

i) for $r\geq i\geq 1$, $F_i$ and $F_{i-1}$ are adjacent chambers
inside $\Bbb A$ ;

ii) $\#{\cal M}(F_B,F'_0)=r.$

\bigskip

Now, if we fix $\tau=(t_r \subset t'_r \supset t_{r-1} \subset
\cdots \subset t'_1 \supset t_0 \subset t'_0)$ a minimal chamber
gallery type, then $t_i=typ (B)=\emptyset$ for $r\geq i\geq 0$
and $t'_j = \{ s_{i_j}\}$ where $s_{i_j}\in S$ is a generator of
$W$ and $t'_0$ is any type. The word $s_{i_r}\cdots s_{i_1}$
built with these reflections, is a reduced decomposition of an
element of $W$.

 More precisely, we have the following
\begin{prop}

i) We keep the notations and $\tau$ is fixed as above. If $F_r$
and $F_0$ are two chambers of $\Bbb A$ such that $d(F_r,F_0) =
r\geq 0$ and $F'_0$ is the face of type $t'_0$ contained in
$F_0$, then $\Gamma (\Bbb A,\tau ,F_r, F'_0) = \set{\g_w}$ is
reduced to a single element.

ii) For any type $t_P$ ($t_P$ is the type of a parabolic subgroup
$P\supset B$), the construction above states a bijection between
minimal chamber gallery types and reduced decompositions of
minimal length coset representatives of the classes in
$W/W_{t_P}$.

\end{prop}

We define an operation on the galleries of $\Gamma (\Bbb A,\tau
,F_r, -)$ by applying to all the faces of a gallery a folding
towards $F_r$. Therefore, all these galleries can be obtained from
$\g_w$, by a finite number of foldings towards $F_r$.


\subsection{Configuration Variety}

Now, we use the combinatorial setting of the building and of the
galleries to describe the Bott-Samelson variety.

Let $P\supset B$ be a parabolic subgroup of $G$ and let $\overline
w\in W/W_{t_P}$. Also, let $w$ denote the element of minimal
length in $\overline w$, and let $w= s_{i_r}\cdots s_{i_1}$ be a
reduced decomposition of $w$, denoted by $\underline i
=(i_r,...,i_1)$.

{\it In the rest of this paper, $w$ will always mean this
element.}

 We fix $F'_0 =F_{\overline w (P)}$, where $\overline w (P)$ stands for
 $\overline w P \overline w^{-1}$. Let

$$\tau_{\underline i}=(t_r \subset t'_r \supset t_{r-1} \subset
\cdots \subset t'_1 \supset t_0 \subset t'_0)$$ be the type of
chamber gallery in $\Bbb A$ between $F_B$ and $F_{\overline w
(P)}$ associated to the reduced decomposition ${\underline i}$,
that is $t_j' = \{s_{i_j}\}$, for $r\geq j \geq 1$ and $t_0' =
t_P$. The type $\tau_{\underline i}$ is therefore a minimal
chamber gallery type.

\begin{defiprop}[{[CC, Part I, $\S$6]}]

The configuration variety
\\$CONF_{\tau_{\underline i}}(G)_B$ is defined as the subvariety of the product

$$\Bbb\Pi (\tau_{\underline i}) = \prod_{r\geq i\geq 0} Par_{t_i}(G)\times Par_{t'_i}(G)$$
which consists of all the configurations of parabolic subgroups
of the form :

$$g =
(B_r \subset P_r \supset B_{r-1} \subset \cdots \subset P_1
\supset B_0 \subset P_0),$$ where $B_r=B$ and
$typ(g)=\tau_{\underline i}$.

\end{defiprop}

As we saw it before, these configurations are the chamber
galleries in $\Delta (G)$, of type $\tau_{\underline i}$ and of
source $F_B$.

To each gallery $g\in\Gamma (\Delta (G), \tau_{\underline i},
F_B,-)$, we associate its target which is a face of type
$t'_0=t_P$, associated to a parabolic subgroup in relative
position $\overline u$ with $B$, for $\overline u\leq \overline w$
in $W/W_{t_P}$. Thus, we obtain a morphism of variety :

$$\pi: CONF_{\tau_{\underline i}}(G)_B \to \overline\Sigma (B, \overline w).$$

Beside, from the decomposition $\underline i =(i_r,\cdots i_1)$,
we can construct the Bott-Samelson variety $\hat\Sigma (B,
\tau_{\underline i})= P_{i_r}\times_B\cdots \times_B (P_{i_1}/B)$,
where for $r\geq j\geq 1$, $P_{i_j} = Bs_{i_j}B\cup B$ is the
unique parabolic subgroup of type $t'_j$ containing $B$ (see for
example [De, $\S$3, Definition 1]). We denote by $[x_r,... ,x_1]$
the $k-$points of this variety.

\begin{prop}[{[CC, Part I, $\S$6]}]

The morphism
$$f: \hat\Sigma (B, \tau_{\underline i})\to CONF_{\tau_{\underline i}}(G)_B$$
defined in the following way, is an isomorphism :

$$f([x_r,\cdots ,x_1])= (B_r \subset P_r \supset B_{r-1} \subset
\cdots \subset P_1 \supset B_0 \subset P_0) $$ where

$B_i = x_r\cdots x_{i+1} B (x_r\cdots x_{i+1})^{-1}$, for
$r-1\geq i\geq 0$,

$P_j = x_r\cdots x_{j+1} P_{i_j} (x_r\cdots x_{j+1})^{-1}$, for
$r-1\geq j\geq 1$,

$B_r=B$ and $P_0 = x_r\cdots x_1 P (x_r\cdots x_1)^{-1}$.

\end{prop}

Thus, we arrive to the well-known

\begin{theo}[{[Ha], [De], [CC, Part I, $\S$6]}]

The couple $\big (CONF_{\tau_{\underline i}}(G)_B,\pi\big )$ is a
smooth equivariant resolution of the Schubert variety
$\overline\Sigma (B, \overline w)$, called the Bott-Samelson
resolution.

\end{theo}

In particular, $\pi$ is a birational proper morphism.

\begin{rem}
Demazure [De] and Hansen [Ha] have first shown independently that
the Bott-Samelson variety is a smooth resolution of the schubert
variety, but Contou Carrere [CC] enunciates this theorem in the
more general setting of the universal Schubert scheme over any
base scheme.
\end{rem}


\section{Combinatorial Tangent Space}

From now on, we will denote the configuration variety as the
Bott-Samelson variety by $\hat\Sigma (B, \tau_{\underline i})$.
And we will call the galleries of $\Gamma (\Bbb A
,\tau_{\underline i}, F_B,-)$, {\it combinatorial galleries}.
Furthermore, we will use freely the two writings of galleries i.e.
$\g = [\g_r,...,\g_1] = (F_B=F_r \supset F'_r \subset F_{r-1}
\supset \cdots \supset F'_1 \subset F_0 \supset F'_0)$. Let us
remark that $\g = [\g_r,...,\g_1]$ is a combinatorial gallery iff
$\g_j = s_{i_j}$ or $\g_j = 1$.


\subsection{The tangent space at $\gamma$ of $\hat\Sigma (B, \tau_{\underline i})$}

Let $Q\supset T$ be a parabolic subgroup of $G$. As $G$ is
``épinglé'', we know that $Lie(Q) = {\cal T}\oplus\big (
\bigoplus_{\alpha\in R_Q} {\cal G}^{\alpha}\big )$, where $R_Q$ is
the set of roots of $(Q,T)$ and where ${\cal G}^{\alpha}=
<X_{\alpha}>$.

If $x$ denote the $k-$point of $Par_{t_Q} (G)$ given by $Q$, then

$$T_x Par_{t_Q} (G)\simeq Lie(G)/Lie(Q) = \bigoplus_{\alpha\in
R\setminus R_Q} {\cal G}^{\alpha}.$$

Let $\tau_{\underline i}=(t_r \subset t'_r \supset t_{r-1} \subset
\cdots \subset t'_1 \supset t_0 \subset t'_0)$ be the type of a
minimal chamber gallery given by a fixed reduced decomposition of
$w\in \overline w$. Let also $\gamma = (F_{B_r} \supset F_{P_r}
\subset F_{B_{r-1}} \supset \cdots \supset F_{P_1} \subset
F_{B_0} \supset F_{P_0})$ denote a gallery of $\Gamma (\Bbb A
,\tau_{\underline i}, F_B,F_{\overline u(P)})$, where $\overline
u\leq\overline w$ in $W/W_{t_P}$ ($P_0=\overline u (P) =
\overline u P \overline u^{-1}$).

From the definition-proposition 1, we deduce that
$T_{\gamma}\hat\Sigma (B, \tau_{\underline i})$ is a vector
subspace of the product :

$$ \Bbb P (\gamma) = \big ({\cal G}/Lie(B_r)\times {\cal
G}/Lie(P_r)\big )\times\cdots \times \big ({\cal G}/Lie(B_0)\times
{\cal G}/Lie(P_0)\big ). $$ Now, we explain, following [CC, Part
II, $\S$5] how to construct a basis of this vector space.

Let $E(\tau_{\underline i}) = \{t_r,t'_r,t_{r-1},\cdots ,t'_1,
t_0, t'_0\}$ be the set of types which appears in
$\tau_{\underline i}$ totally ordered by the reading order (i.e.
$t_r\geq t'_r\geq t_{r-1}\cdots$). Let

$$\Bbb E (\tau_{\underline i}) = \coprod_{\theta\in
E(\tau_{\underline i})} \{\theta\}\times (R\setminus R_{\theta})
\quad {\mathrm {where }} \quad R_{\theta} = \left\{
\begin{array}{rl}
R_{B_i} & {\mathrm {if }} \quad\theta = t_i \\
R_{P_j} & {\mathrm {if }} \quad\theta = t'_j.
\end{array}
\right.
$$

We define an equivalence relation $\sim$ on $\Bbb E
(\tau_{\underline i})$ :

$$(\theta ,\alpha)\sim (\theta' ,\alpha') \Leftrightarrow
\left\{
\begin{array}{l}
 \alpha = \alpha'\\
 \alpha\in (R\setminus R_{\theta})\cap (R\setminus R_{\theta-1})\cap \cdots \cap
 (R\setminus R_{\theta'})\quad (\theta\geq \theta').
\end{array}
\right.
$$

Define $\Bbb B (\gamma) = \Bbb E (\gamma) /\sim$ and $\Bbb B
(\gamma , B) = \{ C\in \Bbb B(\gamma), C\cap \big ( \{t_r\}\times
(R\setminus R_B)\big ) = \emptyset \}$. For $C\in \Bbb B (\gamma ,
B)$, denote $X_C$ the element of $\Bbb P (\gamma)$ with $\theta
-$component equals to $X_{\alpha}$ if $(\theta ,\alpha )\in C$ and
0 otherwise.

\begin{prop}[{[CC, Part II, $\S$5]}]
The set of vectors $\{X_C\}_{C\in \Bbb B (\gamma , B)}$ is a basis
of $T_{\gamma}\hat\Sigma (\tau_{\underline i},B)$.

\end{prop}

\begin{rem}

For $C\in \Bbb B (\gamma , B)$, a vector $X_C$ has the form :

$$ X_C = \big
((0,0),...,(0,0),X_{\alpha},...,X_{\alpha},\underbrace{0,0,...,0}_{\mathrm
{eventually}}\big ). $$

\end{rem}


\subsection{Combinatorial Tangent Space}

Let $x$ be the $k-$point of $\overline\Sigma (B, \overline w)$
given by a parabolic subgroup of type $t_P$ containing $T$. That
is $x= \overline u (P)$ for $\overline u\leq\overline w$ in
$W/W_{t_P}$. Let us denote $\Gamma= \Gamma (\Bbb A
,\tau_{\underline i}, F_B,F_x)$, the finite set of galleries in
$\Bbb A$, of type $\tau_{\underline i}$, of source $F_B$ and of
target $F_x$. These combinatorial galleries are the $T-$fixed
points of the fibre $\hat\Sigma (B,\tau_{\underline
i})_x=\pi^{-1} (x)$ (we shall discuss more precisely this point
in section 4).

\begin{defi}
The combinatorial tangent space, $T_x^c \overline\Sigma (B,
\overline w)$, is the vector subspace of $\overline\Sigma (B,
\overline w)$ generated by the union of the images of tangent
spaces $T_{\gamma}\hat\Sigma (B, \tau_{\underline i})$ for all
$\gamma\in\Gamma$,

$$T_x^c \overline\Sigma (B, \overline w)= \Big <
\bigcup_{\gamma\in\Gamma} T_{\gamma} \pi \big (T_{\gamma}
\hat\Sigma (B, \tau_{\underline i})\big ) \Big >.$$

\end{defi}

This is a vector subspace of $Lie(G)/Lie(\overline u (P)) =
\bigoplus_{\alpha\in R\setminus R_{\overline u (P)}} {\cal
G}^{\alpha}$ and by construction we know a basis of it, composed
by the $X_{\alpha}$ with $\alpha\in R\setminus R_{\overline u
(P)}$, which comes from a $X_C$, such that $X_C = \big
((0,0),...,(0,0),X_{\alpha},...,X_{\alpha}\big )$.

So, $T_x^c \overline\Sigma (B, \overline w)= \bigoplus_{\alpha\in
{\cal R}_x} {\cal G}^{\alpha},$ where ${\cal R}_x\subset
R\setminus R_{\overline u (P)}$ defines the combinatorial tangent
space as a $T-$weighted subspace of $Lie(G)/Lie(\overline u (P))$.
\bigskip

We give here some definitions to understand how a gallery will
contribute to the construction of the combinatorial tangent space.

\bigskip
{\it From now on, we will only consider galleries inside $\Bbb
A$.}

\begin{defi}

A bend $(F_{B'},F_{P'})$ is a gallery of the form $F_{B'}\supset
F_{P'} \subset F_{B'}$, where $F_{B'}$ is a chamber of $\Bbb A$
and $F_{P'}$ is a codimension 1 face of $F_{B'}$ carried by a
wall $M_{\beta'}$.

\end{defi}

\begin{rem}
A generator vector $X_{\beta}$ of $T_x^c \overline\Sigma (B,
\overline w)$ comes from a gallery $\gamma\in\Gamma$ and $\beta$
is associated to a wall $M_{\beta}$ of $\Bbb A$. Two cases arise

1) $\gamma$ crosses $M_{\beta}$ in such a way that $X_{\beta}$
appears in the last component of a $X_C$, $C\in\Bbb B (\gamma ,B)$
(this is the case for all the walls separating $F_B$ from $F_x$);

2) $\gamma$ contains a bend $(F_{B'},F_{P'})$ where $F_{P'}$ is
carried by $M_{\beta}$, again in such a way that $X_{\beta}$
appears in the last component of an $X_C$, $C\in\Bbb B (\gamma
,B)$. Such bends are called {\it generating bends}.

A chamber gallery which contains no bend is minimal, and it
crosses one and only one time each wall which separates its source
from its target.
\end{rem}

\bigskip

As an example, let us suppose that $G$ is the simple group of
type $B_2$ and let us take $P=B$ the "standard" Borel subgroup and
$w = s_{\alpha}s_{\beta}s_{\alpha}$ in the Weyl group which is
the dihedral group generated by $s_{\alpha}$ and $s_{\beta}$. Let
us denote by $P_{\alpha}$ (resp. $ P_{\beta}$) the unique
parabolic subgroup of type $\set{s_{\alpha}}$ (resp.
$\set{s_{\beta}}$) containing $B$.

The following figure shows the graphic representation of the
apartment $\A$ associated to the "standard" torus $T\subset B$ and
of the root system of type $B_2$. We draw the eight combinatorial
galleries of $\Gamma (\Bbb A ,\tau_{\underline i}, F_B,-)$
($\tau_{\underline i}$ is the type associated to the reduced
decomposition $w  = s_{\alpha}s_{\beta}s_{\alpha}$) as curves in
the graphic representation of $\A$. For example $\gamma_0^0$ is
the unique gallery of type $\tau_{\underline i}$ contained in the
chamber $F_B$, so its curve stay inside the chamber while bending
the walls corresponding to the faces $F_{P_{\alpha}}$ and
$F_{P_{\beta}}$. And the curve corresponding to
$\gamma_{\alpha\beta} = [s_{\alpha},s_{\beta},1] = \big
(F_B\supset F_{P_{\alpha}}\subset F_{s_{\alpha}(B)}\supset
F_{s_{\alpha}(P_{\beta})}\subset F_{s_{\alpha}s_{\beta}(B)}\supset
F_{s_{\alpha}s_{\beta}(P_{\alpha})}\subset
F_{s_{\alpha}s_{\beta}(B)}\big )$ follows the faces which appear
in its second writing. And so on for the other galleries. The
bullet at the beginning of each curve stands for the source of
the gallery and the arrow at the end shows the direction and the
target.

\begin{figure}[H]
\centering

\setlength{\epsfxsize}{12cm}

\epsffile{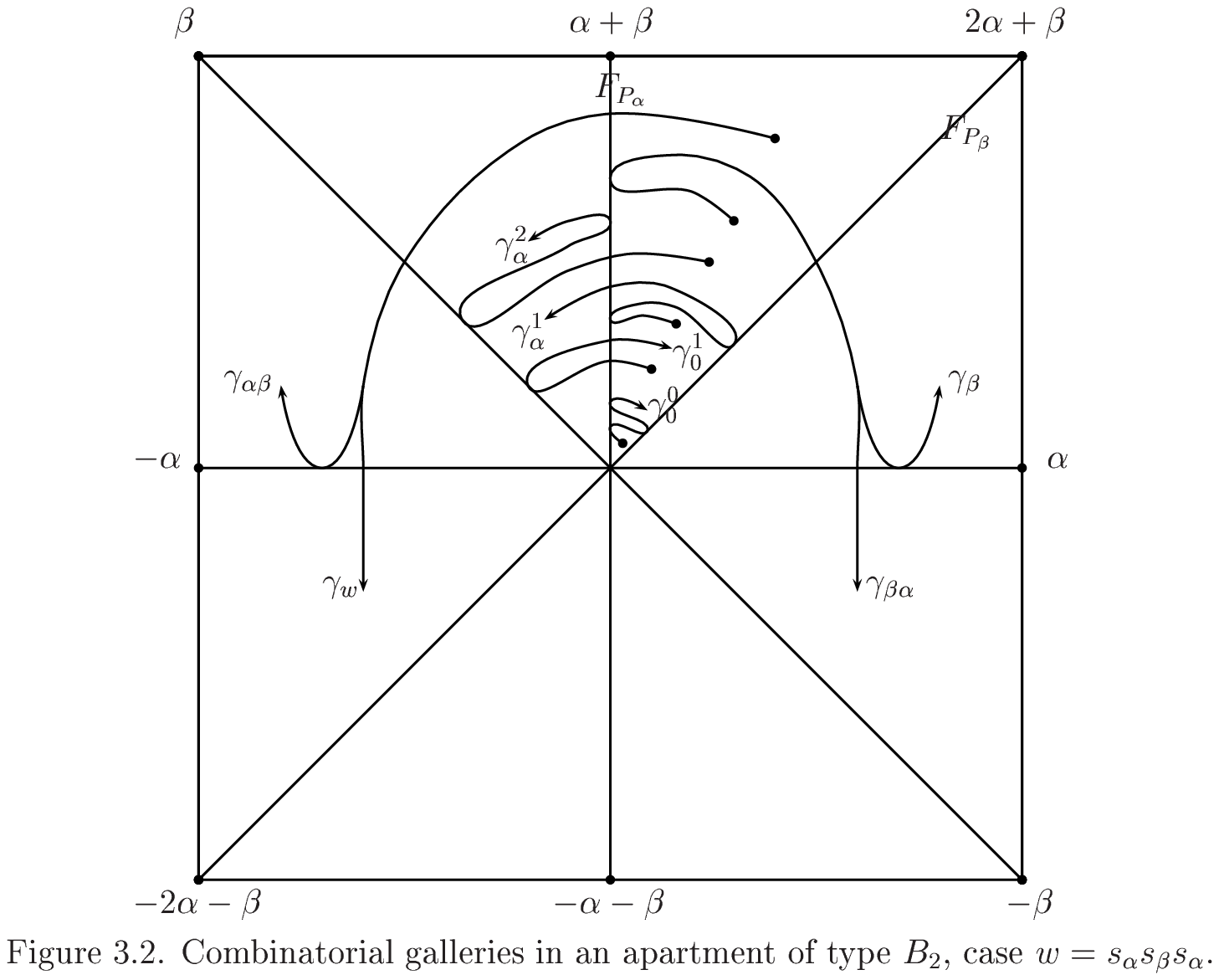}
\end{figure}

As we can see on this picture, if $x = s_{\alpha} B s_{\alpha}$,
the combinatorial tangent space $T^c_x\overline\Sigma(B,w)$ admits
the basis $\set{X_{-\beta}, X_{-2\alpha -\beta},X_{\alpha}}$. And
these generators are all given by generating bends of the two
galleries of $\Gamma (\Bbb A ,\tau_{\underline i}, F_B,F_x) =
\set{\gamma_{\alpha}^1,\gamma_{\alpha}^2}$.

On the other hand, if $x = s_{\alpha}s_{\beta} (B) =
s_{\alpha}s_{\beta} B (s_{\alpha}s_{\beta})^{-1}$, then
$T_x\overline\Sigma(B,w) = T^c_x\overline\Sigma(B,w) =
<X_{\alpha},X_{2\alpha +\beta},X_{-\alpha -\beta}>$. This time,
the two first generators are given by the fact that
$\gamma_{\alpha\beta}$ crosses the walls $M_{\alpha}$ and
$M_{2\alpha +\beta}$ and the third by the generating bend
$(F_{s_{\alpha}s_{\beta}(B)},F_{s_{\alpha}s_{\beta}(P_{\alpha})})$.

\bigskip
It seems that $T_x^c \overline\Sigma (B, \overline w)$ depends on
the choice of a type $\tau_{\underline i}$, i.e. on the choice of
a reduced decomposition of $w$, but it is not the case.

\begin{prop}

The combinatorial tangent space $T_x^c \overline\Sigma (B,
\overline w)$ is independent from the choice of the reduced
decomposition of $w$.

\end{prop}

\noindent {\bf Proof.}\quad Take two decompositions $\underline i$
and $\underline j$ of $w$. The decomposition $\underline j$ can be
obtained from $\underline i$ by the action of an element of the
generalized braid group. This element can be written as a finite
sequence of generalized braid relations. So, it suffices to show
that generators of $T_x^c \overline\Sigma (B, \overline w)$ are
the same if we ``calculate'' them with $\underline i$ or with a
transform of $\underline i$ by a generalized braid relation $r_t$.
Let $\underline j = r_t (\underline i)$, the relation $r_t$
induces a $T-$equivariant isomorphism of variety

$$f_{r_t} : \hat\Sigma (\tau_{\underline i},B)
\stackrel{\simeq}{\longrightarrow} \hat\Sigma (\tau_{\underline
j},B),$$ (hence, $f_{r_t} \big (\Gamma(\tau_{\underline i})\big )
=\Gamma(\tau_{\underline j})$). Let $X_{\beta}$ be a generator of
$T_x^c \overline\Sigma (B, \overline w)$ coming from
$\gamma\in\Gamma (\tau_{\underline i})$. Whether $\gamma$ contains
a minimal chamber gallery $\delta$ which crosses $M_{\beta}$, then
$f_{r_t} (\gamma )$ will contain $f_{r_t} (\delta )$ which still
crosses $M_{\beta}$. Whether $\beta$ is given by a generating
bend $\mu$ of $\gamma$, again $f_{r_t} (\mu )$ will be a
generating bend of $f_{r_t} (\gamma )$ and will still give
$X_{\beta}$ as a generator of $T_x^c \overline\Sigma (B, \overline
w)$. So, in the two cases, proposition holds.

$\hfill\square$
\bigskip

In the following, we will study more closely the set $\Gamma (\A
,\tau_{\underline i} , F_B,F_x)$ of all the combinatorial
galleries above $F_x$. This will allow us to show a pleasant
property of the combinatorial tangent space. In order to do that,
we need some definitions.

\begin{defi}

Let $\gamma = [\g_r,...,\g_1] = (F_B=F_r \supset F'_r \subset
F_{r-1} \supset \cdots \supset F'_1 \subset F_0 \supset F'_0)\in
\Gamma (\A ,\tau_{\underline i} , F_B,F_x)$. A buckle of $\g$ is
a sub-gallery of $\g$ of the shape
$$\mathfrak b = (F'_j \subset
F_{j-1} \supset \cdots \supset F'_{i+1} \subset F_i \supset
F'_i),$$ where $r\geq j>i\geq 0$, such that $F'_j$ and $F'_i$ are
two faces of the same wall $M^{\mathfrak b}$ and such that all the
chambers $F_l$, for $j-1\geq l\geq i$, are on the same side of
$M^{\mathfrak b}$.

\end{defi}

If all the reflections $s_{i_j}$ that compose the gallery of types
$\tau_{\underline i}$ are distinct (including those given by the
type $t_P$) then any gallery $\g$ of $\Gamma (\A
,\tau_{\underline i} , F_B,F_x)$ does not contain any buckle.
Furthermore, any minimal gallery does not contain any buckle
neither. Actually, if all the walls met by a gallery $\g$,
${\mathcal M} (\g) = \set{ M_{\be_r},...,M_{\be_1}, \{ M_{x(\nu
)}\}_{s_{\nu}\in t}}$, are distinct, then $\g$ does not contain
any buckle.

Moreover, a buckle of a combinatorial gallery contains at least
one bend (Definition 7).

\begin{defi}
We say that a buckle $\mathfrak b$ of a gallery $\g$ is maximal
if all the buckles of $\g$ on the wall $M^{\mathfrak b}$ are
contained in $\mathfrak b$ (i.e. are sub-galleries of $\mathfrak
b$).

\end{defi}

Now, we describe the set $\Gamma (\A ,\tau_{\underline i} ,
F_B,F_x)$ using an operation analogous to the one that allowed us
to describe the set $\Gamma (\A ,\tau_{\underline i} , F_B,-)$ in
the section 2.4.

\begin{defi}

A part-folding $\overline\phi$ of a combinatorial gallery
$\g\in\Gamma (\A ,\tau_{\underline i} , F_B,F_x)$ which contains
a buckle $\mathfrak b$ is the following operation. We take the
part of $\g$ determined by $\mathfrak b$ and we replace it by its
image by the reflection $s_{M^{\mathfrak b}}$ associated to the
wall $M^{\mathfrak b}$. We thus obtain a new combinatorial
gallery $\overline\phi(\g)\in\Gamma (\A ,\tau_{\underline i} ,
F_B,F_x)$.

Further, we say that the part-folding $\overline\phi$ is towards
$F_B$ if all the chambers of the new buckle $s_{M^{\mathfrak
b}}(\mathfrak b)$ and $F_B$ are on the same side of $M^{\mathfrak
b}$.

\end{defi}

The set $\Gamma (\A ,\tau_{\underline i} , F_B,F_x)$ is invariant
by part-foldings ($x = \overline u(P)$). Let us denote
\begin{equation}
\g_{u,w}\in \Gamma (\A ,\tau_{\underline i} , F_B,F_x)
\end{equation}
the unique gallery obtained from any other gallery of $\Gamma (\A
,\tau_{\underline i} , F_B,F_x)$ by applying the maximum of
part-foldings towards $F_B$. This gallery corresponds to the
unique way to construct a reduced sub-decomposition of $u$ from
the chosen reduced decomposition of $w$ by taking the reflections
which define $u$ the most possible towards the right.

So, we can recover all the combinatorial galleries above $F_x$ by
applying some part-foldings opposed to $F_B$ (i.e. $F_B$ does not
belong to the image of the associated folding) to $\g_{u,w}$. For
example, if $u = 1$, $\g_{1,w}$ is equal to $\g_0 = [1,...,1]$,
the unique gallery of type $\tau_{\underline i}$ contained in the
chamber $F_B$ and all the galleries having $F_1\subset F_B$ as
target can be obtained from $\g_{1,w}$ thanks to part-foldings
opposed to $F_B$.

\bigskip

The gallery $\g_{u,w}$ carries all the generator vectors of the
combinatorial tangent space $T_x^c \overline\Sigma (B,\overline
w)$.

\begin{prop}

Still keeping the same notations,
$$dimT_x^c \overline\Sigma (B,\overline w) \geq r = dim\overline\Sigma (B,\overline w).$$

\end{prop}

\noindent{\it Proof. } Set $\g_{u,w} = [x_r,...,x_1]\in
P_{i_r}\times_B\cdots \times_B (P_{i_1}/B)$. Recall that the set
of all the walls met by $\g_{u,w}$ is denoted ${\mathcal M}
(\g_{u,w}) = \set{ M_{\be_r},...,M_{\be_1}, \{ M_{u(\nu
)}\}_{s_{\nu}\in t}}$ , where $\be_j = x_r\cdots x_j (\al_{k_j})$
and $u = x_r\cdots x_1$.

Let us suppose that all the walls met by $\g_{u,w}$ are not
distinct (otherwise $\Gamma (\A ,\tau_{\underline i} , F_B,F_x) =
\set{\g_{u,w}}$ and we have $dimT_x^c \overline\Sigma (B,
\overline w) = r$).

On one side, let $\mathcal T (\g_{u,w})\subset\set{M_{\be_r},
...,M_{\be_1}}$ denote the set of walls crossed by $\g_{u,w}$,
i.e. walls which separate $F_B$ from $F_x$. Then $\#\mathcal T
(\g_{u,w}) = l(u)$ and $\g_{u,w}$ never crosses again one of
these walls, otherwise $\g_{u,w}$ would admits a buckle that we
could fold towards $F_B$. Hence, from the remark 4, we obtain
$l(u)$ generator vectors of $T_x^c \overline\Sigma (B,\overline
w)$. Furthermore, the roots associated to the walls which contain
$F_x$ cannot give any generator of the combinatorial tangent
space. Hence, only the $m = l(w) - l(u)$ walls of $\mathcal M
(\g_{u,w})\setminus \mathcal T (\g_{u,w}) =
\set{M_{\be_{i_m}},...,M_{\be_{i_1}}}$ may provide some new
generators.

On the other side, the hypothesis done (i.e. all the walls met by
$\g_{u,w}$ are not all distinct), implies that there exists only
$p$ walls, $\mathcal D(\g_{u,w}) =
\set{M_{\be_{i'_p}},...,M_{\be_{i'_1}}}\subset
\set{M_{\be_{i_m}},...,M_{\be_{i_1}}}$ with $m>p\geq 1$, all
distinct and also distinct from the walls crossed by $\g_{u,w}$
on which this gallery admits a generating bend. So, they give $p$
more generator vectors of $T_x^c \overline\Sigma (B,\overline w)$
(Remark 4).

\bigskip

The rest of the walls, $\mathcal R (\g_{u,w}) =
\set{M_{\be_{i_m}},...,M_{\be_{i_1}}} \setminus
\set{M_{\be_{i'_p}},...,M_{\be_{i'_1}}}$ are then equal to some
crossed walls or some walls containing $F_x$ or some walls of
$\mathcal D(\g_{u,w})$. Hence, $\#\mathcal R (\g_{u,w}) = m-p = t
+ f + c$, where $t$ is the number of walls of $\mathcal R
(\g_{u,w})$ equal to crossed walls, $f$ those equal to walls
containing $F_x$ and $c$ those equal to walls of $\mathcal
D(\g_{u,w})$. Moreover, $t = t_1 + \cdots +t_s$, $f = f_1+\cdots
+f_e$ and $c = c_1 +\cdots + c_b$, where the $t_j$'s, $f_j$'s and
$c_j$'s are respectively the numbers of walls of $\mathcal R
(\g_{u,w})$ equal to a same wall of $\mathcal T (\g_{u,w})$, of
$\{ M_{u(\nu )}\}_{s_{\nu}\in t}$ and of $\mathcal D(\g_{u,w})$.

\bigskip
By the way, the gallery $\g_{u,w}$ admits $s$ maximal buckles on
crossed walls, $e$ maximal buckles on walls containing $F_x$ and
$b$ maximal buckles on walls of $\mathcal D(\g_{u,w})$ (cf.
Definition 9). And each of these buckles provides, by applying
part-foldings opposed to $F_B$, at least as much galleries
containing a generating bend as the number $t_j$, $f_j$ or $c_j$
which corresponds to the maximal buckle we started with.

Indeed, let $\mathfrak b$ be one of these buckle, $\mathfrak b$
owns at least as much bends as the number $t_j$, $f_j$ or $c_j$
which corresponds to it. Let us apply to $\g_{u,w}$ the
part-folding $\overline\phi$ fixed by $\mathfrak b$ (Definition
10), we obtain a new gallery $\overline\phi(\g_{u,w})$ which owns
a buckle, image of $\mathfrak b$. Let $(C',F')$ be a bend of this
new buckle, then two cases may arise :

a) $(C',F')$ is a generating bend on a wall that $\g_{u,w}$ did
not meet ;

b) the gallery $\overline\phi(\g_{u,w})$ crosses the wall
containing $F'$, that means that it has a buckle on this wall,
then we can apply a part-folding to it and repeat the discussion.

At the end, each maximal buckle will give at least as much
generator vectors of $T_x^c \overline\Sigma (B,\overline w)$ as
the number $t_j$, $f_j$ or $c_j$ corresponding to it.

\bigskip
Hence, if we recount the number of generators which may arise
from $\g_{u,w}$, we see that $dimT_x^c \overline\Sigma (\overline
w) \geq r = l(w)$.

$\hfill\square$

\begin{rem}
All the combinatorial galleries above $F_x$ are obtained by
part-foldings opposed to $F_B$ from $\g_{u,w}$. So, all the
generators of $T_x^c \overline\Sigma (B, \overline w)$ will be
given using the process shown in the previous proof. Hence, the
$T-$weights of this space, $\mathcal R_x$, will come from the
walls we consider above.

\end{rem}

\bigskip

To close this section, we mention an important result of C.
Contou Carrere.

\begin{theo}[{[CC, Part II, $\S$10]}]

Let $G$ be a $k-$semi-simple group of type $A$. For any parabolic
subgroup $P\supset B$ and for all $T-$fixed point $x\in
\overline\Sigma (B, \overline w)$, the combinatorial tangent
space $T_x^c \overline\Sigma (B, \overline w)$ is equal to the
Zariski tangent space $T_x \overline\Sigma (B, \overline w)$.

\end{theo}

Thus, we can describe the singular locus of any Schubert variety
in case of a group of type $A$.


\section{Rational Smoothness of Schubert Varieties}

First of all, we recall one of the definitions of rational
smoothness [CP]. A projective variety $Y$ of dimension $r$ is
rationally smooth at $y\in Y$ if there exists an neighborhood $U$
of $y$ such that the $l-$adic cohomology $H^i_{<y>}(U) = 0$ if
$i\not = 2r$ and is one dimensional if $i=2r$. The variety $Y$ is
rationally smooth if it is rationally smooth at every point.

A smooth point is rationally smooth, but there are rationally
smooth Schubert varieties which are not smooth. This is the case
for our example of the figure 3.2, see [C2, $\S$7] for other
examples.

\subsection{Combinatorial Tangent Space and Rational Smoothness}

\begin{prop}

Let $\overline w\in W/W_t$. For all $\overline u\leq \overline
w$, we set $x=\overline u (P)$, then
$$dimT^c_x\overline\Sigma(B,\overline w) = \#\set{\be\in R_+,\
\overline {s_{\be}u}\leq \overline w \ \hbox{ and } \overline
{s_{\be}u}\not = \overline u}.$$

\end{prop}

\noindent{\it Proof. } First, let $\be\in R_+$ such that
$\overline {s_{\be}u}\leq \overline w$ and $\overline
{s_{\be}u}\not = \overline u$, we denote by $M_{\be}$ the wall
associated to this root. Two cases may appear :

1) $M_{\be}$ separates $C$ from $F_x$, the face associated to
$\overline u(P)$. Then $X_{\be}$ is a generator vector of
$T^c_x\overline \Sigma(B,\overline w)$.

2) $C$ and $F_x$ are on the same side of $M_{\be}$. Thus, the
gallery $\g_{s_{\be}u,w}$ (cf. (1) of the previous section for the
definition) crosses this wall. Let us consider the gallery
$\g\in\Gamma (\A ,\tau_{\underline i} , F_B,F_x)$ image of the
folding towards $F_B$ associated to $M_{\be}$ of the gallery
$\g_{s_{\be}u,w}$. Then $\g$ contains a generating bend on this
wall and that gives a generator vector $X_{-\be}$ of the
combinatorial tangent space $T^c_x\overline \Sigma(B,\overline
w)$.

Moreover, the face $F_x$ cannot belong to the wall $M_{\be}$
since this is equivalent to the fact that $\overline {s_{\be}u} =
\overline u$.

\bigskip
Second, let $M_{\al}$ be a wall that provides a generator vector
of the combinatorial tangent space (i.e. $X_{\al}$ or $X_{-\al}$
belongs to $T^c_x\overline \Sigma(B,\overline w)$). Let
$\g\in\Gamma (\A ,\tau_{\underline i} , F_B,F_x)$ denote a gallery
which gives this generator. Then the gallery obtained from $\g$ by
keeping $\g$ until the wall $M_{\al}$ and by following it with
$s_{\al}(\g)$ is a combinatorial gallery that admits $F_{\overline
{s_{\al}u}(P)}$ as a target. Hence, $\overline {s_{\al}u}\leq w$
and $\overline {s_{\al}u} \not = \overline u$ since a wall which
gives a generator vector of $T^c_x\overline\Sigma(B,\overline w)$
does not contain $F_x$.

\bigskip
These two constructions are inverse of each other, so they
establish a bijection between the set of generator vectors of
$T^c_x\overline\Sigma(B,\overline w)$ and $\set{\be\in R_+,\
\overline {s_{\be}u}\leq \overline w\ \hbox{ and } \overline
{s_{\be}u}\not = \overline u}$, hence the proposition is proved.

$\hfill\square$

\bigskip

Let $\overline u\leq \overline w$ inside $W/W_t$. Let us
consider, following Carrell and Peterson ([CP]), the property
$\mathbf P(\overline u,\overline w)$ defined by :
$$\mathbf P(\overline
u,\overline w) \Leftrightarrow \forall\overline v,\ \overline
u\leq\overline v\leq\overline w,\ \#\set{\be\in R_+,\ \overline
u<\overline {s_{\be}v}\leq \overline w\ \hbox{ and } \overline
{s_{\be}v}\not = \overline v} = l(\overline w) - l(\overline u),$$
where the lengths are taken on the elements of minimal length of
each class.

As it is remarked in [C1], this property is equivalent to the
following, for all $\overline u\leq\overline v\leq\overline w,\
\#\set{\be\in R_+,\ \overline {s_{\be}v}\leq \overline w\ \hbox{
and } \overline {s_{\be}v}\not = \overline v} = l(\overline w)$.

\bigskip
Now, let us suppose that $P = B$, i.e. we consider the Schubert
varieties inside the variety of Borel subgroups. Then the last
conditions, $\overline {s_{\be}v}\not = \overline v$, in the above
properties, are superfluous.

Let $u,w$ be two elements of the weyl group such that $u\leq w$.
We denote by $P_{u,w}$ the Kashdan-Lusztig polynomial (cf. [KL],
[D90]). We recall the following result due to Carrell and
Peterson.

\begin{theo}[{[CP, Theorem C]}]

Let $u\leq w$ inside the Weyl group $W$, then the property
$\mathbf P(u,w)$ is equivalent to $P_{u,w} = 1$.

\end{theo}

But, in our case, the rational smoothness of a Schubert variety
$\overline\Sigma(B,w)$ at the point $u(B)$ is characterized by
the fact that $P_{u,w} = 1$ (see for exemple [KL] in the case the
ground field has positive characteristic and [CP] for any
algebraically closed field). Thus, as a corollary of what we said
from the begining of this section (in particular, from
Proposition 6), we get :

\begin{theo}

A Schubert variety $\overline\Sigma(B,w)$ is rationally smooth at
$u(B)$, for $u\leq w$, if and only if for all $u\leq v\leq w$, the
dimension of the combinatorial tangent space $T^c_{v(B)}
\overline\Sigma(B, w)$ is equal to $r = l(w)$.

\end{theo}

Furthermore, combining Theorem 2 of section 3.2 and Theorem 4, we
have the following result, first obtained by Deodhar [D85] :

\begin{coro}
If $G$ is of type $A$, then smoothness and rational smoothness for
Schubert varieties in $G/B$ are equivalent.

\end{coro}

\begin{rem}
This characterization seems to go through in the parabolic case,
that is for any $P\supset B$. But, in this case, it has to be
taken into account the parabolic Kashdan-Lusztig polynomials to
characterize the rational smoothness of Schubert varieties in
$G/P$ (cf. [D91, $\S$6]).

\end{rem}

\subsection{Combinatorial Tangent Space and $T-$invariant Curves}

In this section, we keep the notations as in section 3 and we
follow the definitions of Carrell and Peterson [CP] about
$T-$invariant curves.

For any projective $T-$variety $Y$, a $T-$invariant curve in $Y$
is the closure of a one dimensional orbit of $T$. Let us denote
by $E(Y)$ the set of all the $T-$invariant curves in $Y$ and if
$x$ is a $T-$fixed point of $Y$, $E(Y,x)$ denotes the set of all
the $T-$invariant curves containing $x$.

First, we recall here some well-known properties of the
$T-$invariant curves in
$$\overline\Sigma (B, \overline w)\hookrightarrow Par_{t_P} (G) \simeq G/P.$$

For a positive root $\alpha\in R^+$, we denote by $U_{\alpha}$ the
unipotent subgroup of rank 1 of $G$ associated to $\alpha$. And
$Z_{\alpha} = <U_{\alpha}, U_{-\alpha}> \subset G,$ denotes the
copy of $Sl_2$ associated to $\alpha$.

Every $T-$invariant curve ${\mathcal C}$ in $Par_{t_P}(G)$ has
the form ${\mathcal C}= Z_{\alpha}\overline w (P)$ for some
$\overline w\in W/W_{t_P}$ and $\alpha\in R^+$ [CP, Theorem F].

Let $\overline u\in W/W_{t_P}$ such that $\overline u\leq\overline
w$ for the relative Chevalley-Bruhat order (defined on the set of
minimal length coset representatives for $W/W_{t_P}$). Set
$x=\overline u (P)$, the $T-$fixed point corresponding to
$\overline u$.

Carrell and Peterson proved that $Z_{\alpha} \overline u (P)\in
E\big (\overline\Sigma (\overline w,B) , x \big )$ if and only if
$\overline u\leq \overline w$ and $\overline {s_{\alpha} u} \leq
\overline w$ (and $\overline {s_{\alpha} u}\not =\overline u$)
[CP, Theorem F].

\bigskip

Thus, from the proposition 6 (more precisely from the proof of
this proposition), we obtain,

\begin{prop}
To each generator vector of the combinatorial tangent space
$T^c_x\overline\Sigma(B, \overline w)$ corresponds a
$T-$invariant curve in $\overline\Sigma (\overline w)$ which
contains $x$ and this space is spanned by the tangent vectors to
each of these curves. Hence, $\#E\big (\overline\Sigma (\overline
w) , x \big ) = dim T^c_x\overline\Sigma(B, \overline w)$.

\end{prop}

\bigskip
Now, we go to the $T-$invariant curves in the Bott-Samelson
variety $\hat\Sigma (B, \tau_{\underline i})$.

The set of $T-$fixed points of this variety $\hat\Sigma (B,
\tau_{\underline i})^T$ is finite and is equal to the finite set
$\Gamma (\Bbb A, \tau_{\underline i}, F_B, -)$ of all the
combinatorial galleries of $ \hat\Sigma (B, \tau_{\underline
i})$. Hence, there are at least $r$ $T-$invariant curves passing
through each combinatorial gallery [CP, $\S$2].

\begin{defi}
A $T-$invariant curve $\mathfrak C$ of the Bott-Samelson variety
will be called a combinatorial $T-$invariant curve if it contains
two combinatorial galleries $\gamma =[\g_r,...,\g_1]$ and $\delta
= [\delta_r,...,\delta_1]$ for which there exists a unique $j_0\in
\{ r,..., 1\}$ such that for all $j\not = j_0$, $\g_j = \delta_j$
and $\g_{j_0}= s_{\al_{k_{j_0}}} \delta_{j_0}$.

\end{defi}

Such a curve $\mathfrak C$ maps by $\pi$ on a $T-$invariant curve
in the Scubert variety, more precisely $\pi(\mathfrak C) =
Z_{\beta} \pi (\gamma ) = Z_{\beta} \pi (\delta )$, where $\beta =
\g_r\cdots \g_{j_0 +1} (\alpha_{k_{j_0}})$. However, there exists
some $T-$invariant curves in the Bott-Samelson variety that are
not of this kind. This is due to the fact that the action of the
torus $T$ on $ \hat\Sigma (B, \tau_{\underline i})$ is not
"special" in the sense of Carrell and Peterson [CP, $\S$2].

So, we use the notation ${\mathfrak C} = (\gamma\frown\delta )$ to
describe such a $T-$invariant curve, which contains $\gamma$ and
$\delta$ as two distinct $T-$fixed points.

Moreover, some of the galleries of $\Gamma (\Bbb A
,\tau_{\underline i}, F_B , -)$ can be obtained from the gallery
$\gamma_{w}\in\Gamma (\Bbb A ,\tau_{\underline i}, F_B
,F_{\overline w (P)})$ by a finite number of foldings towards
$F_B$ (cf. 2.4, Proposition 1). Hence, if ${\mathfrak C} =
(\gamma\frown\delta )\in E\big (\hat\Sigma (B, \tau_{\underline
i})\big )$, where $\delta$ is obtained from $\gamma$ by a folding
$\phi_{\alpha}$ towards $F_B$ ($\alpha\in R^+$), we will denote
it by ${\frak C}_{\alpha} = (\gamma\stackrel
{\phi_{\alpha}}{\frown}\delta ).$ Furthermore, $\pi ({\mathfrak
C}_{\alpha}) = Z_{\alpha} \pi (\gamma ) = Z_{\alpha} \pi (\delta
)$.

\bigskip

We now define the finite graph ${\cal BS}^c_T$ of the
combinatorial $T-$invariant curves in $\hat\Sigma (B,
\tau_{\underline i})$.

The set of vertices ${\cal BS}^c_T$ is defined as the set
$\hat\Sigma (B, \tau_{\underline i})^T$ of all $T-$fixed points
in the Bott-Samelson variety, i.e. the combinatorial galleries
$\Gamma (\Bbb A ,\tau_{\underline i}, F_B ,-)$.

And the set of edges ${\cal BS}_T^c$ is defined as the set of
combinatorial curves in $E\big (\hat\Sigma (B, \tau_{\underline
i})\big )$. For $\gamma$ and $\delta$ in $\Gamma (\Bbb A
,\tau_{\underline i}, F_B ,-)$, there is an edge between them if
${\frak C} = (\gamma\frown\delta )$ is a combinatorial
$T-$invariant curve in $\hat\Sigma (B, \tau_{\underline i})$. If
this curve is associated to a folding $\phi_{\alpha}$ towards
$F_B$ ($\al>0$), we label the corresponding edge of ${\cal
BS}_T^c$ by $\phi_{\alpha}$.

Moreover, we form into rows the vertices of the graph (i.e. the
combinatorial galleries) depending on the distance that separate
the chamber $F_B$ from the target of the gallery. Thus the
gallery $\g_w$ will be the unique vertex at the bottom row of the
graph and $\g_{1,w}$ will belong to the top row. As this distance
is exactly the length of the element of the Weyl group that gives
the target face of the gallery, we have a surjective morphism of
graphs from ${\cal BS}_T^c$ to the Bruhat subgraph of $(W,S)$
associated to $w$, hence we say that ${\cal BS}_T^c$ is above the
Bruhat graph of $w$.

Furthermore, from the previous definition, exactly $r=
dim\hat\Sigma (B, \tau_{\underline i})$ edges meet at each vertex.

\bigskip
To see what such a graph looks like, let us take the example of
$\S3.2$ up again. We recall that the group $G$ is of type $B_2$
and we fix $P=B$ and $w = s_{\alpha}s_{\beta}s_{\alpha}$. The set
of combinatorial galleries $\Gamma (\Bbb A ,\tau_{\underline i},
F_B ,-)$ is represented in the figure 3.2. And from their
expressions in terms of generators of the Weyl group as points of
$P_{\alpha}\times_B P_{\beta}\times_B P_{\alpha}/B$, it is easy
to describe the combinatorial $T-$invariant curves between them.
So in this case the graph ${\cal BS}_T$ has the following form.

\begin{figure}[H]
\centering

\setlength{\epsfxsize}{12cm}

\epsffile{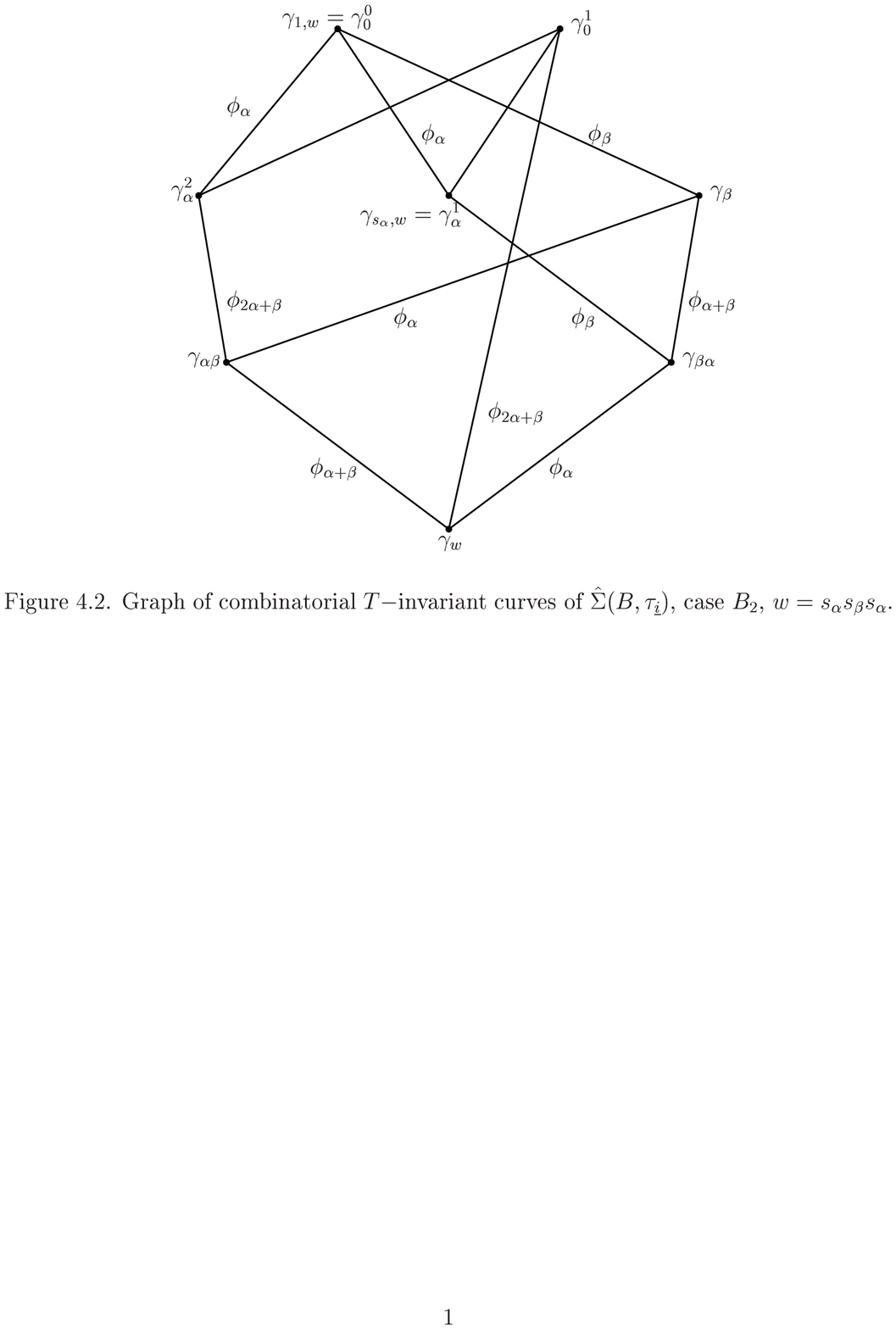}
\end{figure}

As we can see on this example, some of the $T-$invariant curves
of $\hat\Sigma (B, \tau_{\underline i})$ are not related to a
folding. For the others, the folding can be easily read on the
figure 3.2.

\bigskip
Now, we return to the general case and we describe a way to read
off the generators of the combinatorial tangent space from the
graph ${\cal BS}_T^c$.

First, let $x=\overline u (P)$ be a $T-$fixed point of
$\overline\Sigma (B, \overline w)$, i.e. a vertex of the Bruhat
graph of $w$.

Let us consider all the vertices above $x$, that is, all the
combinatorial galleries whose target is $F_x$. We denote by
$\mathcal R_x^+$ (resp. $\mathcal R_x^-$) the set of all the
distinct roots that, eventually, label the edges leaving to the
top (resp. arriving from the bottom to) these vertices. Then, we
deduce from what has been said from the begining of this section
the following proposition.

\begin{prop}

The set ${\mathcal R}_x$ of all the $T-$weights of the
combinatorial tangent space $T^c_x\overline\Sigma(B,\overline w)$
(cf. Section 3.2) actually equals $\mathcal R_x^+ \amalg
(-\mathcal R_x^-)$.

\end{prop}

This proposition shows up a new way to calculate the generator
vectors of the combinatorial tangent space. For instance, if we
go back to the previous example, it is easy to see on the graph
that $T^c_{s_{\al}(B)}\overline\Sigma(B,w)$ admits the basis
$\set{X_{-\beta}, X_{-2\alpha -\beta},X_{\alpha}}$ and if $x =
s_{\alpha}s_{\beta} (B) = s_{\alpha}s_{\beta} B
(s_{\alpha}s_{\beta})^{-1}$, then $ T^c_x\overline\Sigma(B,w) =
<X_{\alpha},X_{2\alpha +\beta},X_{-\alpha -\beta}>$. Thus, we
recover the results of section 3.2.


\vskip 2cm

\section{References}

\vskip 5mm
\begin{description}

\item[{\bf [BOU]}] {\sc N. Bourbaki}, {\em Groupes et algèbre
de Lie, Chapitres 4,5 et 6}, Hermann, Paris (1968).

\item[{\bf [C1]}] {\sc J. B. Carrell}, On the smooth points of
a Schubert variety, {\em Can. Math. Soc. Conf. Proc.}, vol {\bf
16}, (1995), 15-33.

\item[{\bf [C2]}] {\sc J. B. Carrell}, The span of the tangent
cone of a Schubert variety, {\em Algebraic Groups and Lie Groups,
Soc. Lectures Series 9, Cambridge Univ. Press},  {\bf } (1997),
51-60.

\item[{\bf [CP]}] {\sc J. B. Carrell}, The Bruhat graph of a
Coxeter group, a conjecture of Deodhar, and rational smoothness of
Schubert varieties, {\em Proc. Symp. in Pure Math.},  {\bf 56},
(1994), Part I, 53-61.

\item[{\bf [CC]}] {\sc C. Contou Carrere}, {\em Géométrie des
groupes semi-simples, résolutions équivariantes et lieu singulier
de leurs variétés de Schubert}, {\em Thèse d'état}, (1983),
Université Montpellier II, published partly as, Le Lieu singulier
des variétés de Schubert, {\em Adv. Math.}, {\bf 71}, (1988),
186-221.

\item[{\bf [De]}] {\sc M. Demazure}, Désingularisation des
variétés de Schubert, {\em Ann. Sci. Ecole Norm. Sup. (4)}, {\bf
7}, (1974), 53-88.

\item[{\bf [D91]}] {\sc V. Deodhar}, A brief survey of Kashdan-Lusztig
theory and related topics, {\em Proc. Symp. in Pure Math.}, {\bf
56}, (1991), Part I, 105-124.

\item[{\bf [D90]}] {\sc V. Deodhar}, A combinatorial setting for question
in Kashdan-Lusztig theory, {\em Geometriae Dedicata}, {\bf 36},
(1990), 95-119.

\item[{\bf [D85]}] {\sc V. Deodhar}, Local Poincaré duality and
non singularity of Schubert varieties, {\em Comm. in Algebra},
{\bf 13(6)}, (1985), 1379-1388.

\item[{\bf [G]}] {\sc  S. Gaussent}, Etude de la résolution de Bott-Samelson
d'une variété de Schubert, en vue d'un critère valuatif de
lissité, Ph. D. Thesis, Université Montpellier 2, in preparation,
to be defended in january 2001.

\item[{\bf [H]}] {\sc H. C. Hansen}, On cycles of flag
manifolds, {\em Math. Scand.}, {\bf 33}, (1973), 269-274.

\item[{\bf [KL]}] {\sc D. Kazhdan and G. Lusztig},
Representation of Coxeter groups and Hecke algebras, {\em Invent.
Math.},  {\bf 53}, (1979), 165-184.

\item[{\bf [Ma]}] {\sc O. Mathieu},
Formules de caractères pour les algèbres de Kac-Moody générales,
{\em Astérisque}, {\bf 159-160}, (1988).

\item[{\bf [SGA 3]}] {\sc A. Grothendieck and M. Demazure}, {\em
Schémas en groupes I, II, III, Lectures notes in Math. 151, 152,
153}, Springer-Verlag, Heidelberg (1970).

\item[{\bf [T]}] {\sc J. Tits}, {\em Buildings of spherical
type and finite $BN-$pairs, Lecture notes in Math. 386},
Springer-Verlag, Heidelberg (1974).

\end{description}

\end{document}